\documentclass[12pt]{article}

\usepackage[ansinew]{inputenc}
\usepackage{pst-all,epsfig}
\usepackage{amsmath}
\usepackage{amsfonts}
\usepackage{amssymb}

\newtheorem{claim}{Claim}
\newtheorem{lemma}{Lemma}
\newtheorem{cor}{Corollary}

\newtheorem{theorem}[lemma]{Theorem}

\newcommand{\tri}{\triangle}
\newcommand{\al}{\alpha}

\newcommand{\si}{\sigma}
\newcommand{\la}{\lambda}
\newcommand{\eps}{\varepsilon}

\newcommand{\R}{\mathbb{R}}

\newcommand{\ver}{\textrm{vert}}
\newcommand{\conv}{\textrm{conv}}
\newcommand{\intt}{\textrm{int }}

\newcommand{\fin}{\hfill $\Box$}

\begin{document}

\title {Helly type theorems for the sum of vectors in a normed plane}
\author{Imre B\'ar\'any and Jes\'us Jer\'onimo-Castro}

\maketitle

\begin{abstract}
The main results here are two Helly type theorems for the sum
of (at most) unit vectors in a normed plane. Also, we give a new
characterization of centrally symmetric convex sets in the plane.
\end{abstract}


\section{Main results}
This paper is about the sum of vectors in a normed
plane. We fix a norm $\|.\|$ in $\R^2$ whose unit ball is $B$; so
$B$ is an $0$-symmetric convex body. There are some interesting results about sums
of unit vectors in normed planes. For instance, it is proved in
\cite{bar} that for every subset $V=\{v_1,\ldots ,v_n\}\subset B$
of unit vectors, with $n$ an odd number, we may choose numbers
$\epsilon_1, \epsilon_2, \ldots, \epsilon_n$ from $\{1,-1\}$ such
that $\|\sum\limits_{v_i\in V} \epsilon_i v_i \|\leq 1.$ This time we
are interested in unit vectors whose sum has length at least 1.

We write $u\cdot v$ for the usual
scalar product of $u,v\in \R^2$ and $[n]$ for the set
$\{1,2,\ldots,n\}$. Here comes our first result.

\begin{theorem}{\label{semiespacio}}Assume $n\ge 3$ is an odd integer
and $V=\{v_1, v_2, \ldots , v_n\} \subset \R^2$ is a set unit vectors. If $u\cdot
v_i \ge 0$ for every $i\in [n]$ with a suitable non-zero vector $u \in \R^2$,
then
$$\|v_1 + v_2 + . . . + v_n\| \ge 1.$$\end{theorem}

Here and in what follows we can assume that $V$ is a multiset, that is, $v_i=v_j$ can happen even
if $i\ne j$. Perhaps one should think of $V$ as a sequence of $n$ vectors from $\R^2$.

Our main results are two unusual Helly type theorems whose proof uses Theorem~\ref{semiespacio}.
For information about Helly type results the reader may consult \cite{DGK}.

\begin{theorem}{\label{hellytype}}Assume $n\ge 3$ is an odd integer
and $V=\{v_1, v_2, \ldots , v_n\} \subset \R^2$ is a set unit vectors. If the sum
of any three of them has norm at least 1, then
$$\|v_1 + v_2 + . . . + v_n\| \ge 1.$$\end{theorem}

\begin{theorem}\label{helly} Assume $n\ge 3$ is an odd integer
and $V=\{v_1, v_2, \ldots , v_n\} \subset B$. If the sum of any
three elements of $V$ has norm larger than 1, then
$$\|v_1 + v_2 + . . . + v_n\| > 1.$$\end{theorem}

To our surprise Theorem~\ref{helly} fails in the following form: If $V \subset B$, $|V|$ is odd, and
the sum of any three of its elements has norm at least 1, then
$\|v_1 + v_2 + . . . + v_n\| \ge 1.$ The example is with the max norm and the vectors are $v_1=(1,1)$, $v_2=(-1,1)$,
and $v_3=v_4=v_5=(0,-1/2)$. This is also an example showing that Theorem~\ref{hellytype} does not hold
if we require $V \subset B$ instead of $\|v_i\|=1$ for all $i$.

Note that in these theorems $n$ has to be odd. Indeed, let $w_1$
and $w_2$ two almost antipodal unit vectors with $\|w_1+w_2\|$
very small, say $\eps$, set $n=2k$, $v_1=\ldots=v_k=w_1$ and
$v_{k+1}=\ldots=v_n=w_2$. The conditions of our three Theorems are
satisfied (except that $n$ is even now) but $\|v_1 + v_2 + . . . +
v_n\| =\eps n$, as small as you wish if you choose $\eps$ small
enough.

For simpler writing let $[n] \choose k$ denote the set of all $k$-element subsets of $[n]$, and given $S \in {[n] \choose k}$ define
\[
\si (S,V)=\sum_{i \in S}v_i,
\]
and we call it a $k$-sum of $V$. Note that $\si(\emptyset,V)=0$ by definition. Theorem~\ref{helly} has the following immediate

\begin{cor}Assume $n\ge 5$ is an integer, $V=\{v_1, v_2, \ldots , v_n\} \subset B$, $k\in [n]$ is odd and $k >3$. If every 3-sum of $V$ is outside $B$, then so is every $k$-sum of $V$.
\end{cor}

Theorems~\ref{semiespacio} and \ref{hellytype} have similar corollaries and the interested reader will have no difficulty stating or proving them.

We close this section with a neat {\bf proof} of Theorem~\ref{semiespacio} for the case of
the Euclidean norm. The method (unpublished) is due to Boris Ginzburg who used it for the Euclidean case of Theorem 1 from \cite{bar}.

We may assume w.l.o.g. that $u=(0,1)$. The proof is in fact an
algorithm that produces a sequence $V=V_0,V_1,\ldots,V_n$ of sets
of $n$ unit vectors, satisfying $u\cdot v\ge 0$ for all $v \in
V_i$, $i\in [n]$ so that the norm of $s_i=\sum_{v \in V_i}v$
decreases as $i$ increases and $\|s_n\|\ge 1$. Call an element $v
\in V_i$ {\sl fixed} if it equals $(1,0)$ or $(-1,0)$, and let
$F_i$ be the set of fixed elements in $V_i$, and let
$M_i=V_i\setminus F_i$ the set of {\sl moving} elements in $V_i$.

At the start  $V=V_0=M_0$ and $F_0=\emptyset$. Assume $V_i$ has
been constructed, and set $f_i=\sum_{v\in F_i}v$ and $m_i=\sum_{v
\in M_i}v$. One can rotate the vector $m_i$ so that $\|f_i+m_i\|$
decreases during the rotation (because of the cosine theorem). We
rotate $m_i$ in this direction, together with all vectors in $M_i$
as long as one of its elements, say $v^*$, reaches $(1,0)$ or
$(-1,0)$. Let $M_i^*$ be this rotated copy of $M_i$. Define
$M_{i+1}=M_i^*\setminus \{v^*\}$ and $F_{i+1}=F_i\cup\{v^*\}$. We
indeed have $\|s_i\|\ge \|s_{i+1}\|$. By construction $V_n=F_n$,
$M_n=\emptyset$ and $\|f_n\|$ is an odd integer so
$\|s_n\|=\|f_n\|\ge 1$.\fin

\section{Proof of Theorem~\ref{semiespacio}}

\noindent {\bf Proof.} We assume again that $u=(0,1)$. Let
$n=2k-1$ and let $v_1,\ldots ,v_{2k-1}$ be our unit vectors in
clockwise order on the boundary of $B$ in the upper halfplane. Let
$w_1$ and $w_2$ be two unit vectors on the horizontal line through
$0$ with $w_1$ to the left of the origin $0$. The tangent line $L$
to $B$ at $v_k$ bounds the half-plane $H$, the one not containing
the origin.  Set $s=v_1+ \ldots + v_{2k-1}$.

Let $\ell$ be the line through $0$ and $v_k$. For $v \in \R^2$ let $v'$ be the signed length of its projection in direction $L$ onto $\ell$, that is,
$v'$ is positive if $v'$ has the same direction as $v_k$ and negative otherwise.
Since the projection of the sum of vectors is equal to the sum of their
projections, it suffices to prove that $$v'_1+v'_2+ \ldots +v'_{2k-1}\ge
1$$ as this implies $s \in H$ and so $\|s\| \ge 1$. We have that $v'_k=\|v_k\|=1$ and
$$v'_1+\ldots +v'_{k-1}\geq (k-1)w'_1$$
$$v'_{k+1}+\ldots +v'_{2k-1}\geq (k-1)w'_2.$$
As $w'_1+w'_2=0$, the proof is now complete. \fin

\noindent {\bf Remark 1.} Using this proof the case of equality
can be characterized but the conditions are clumsy. The case when
the boundary of $B$ contains no line segment is simple: equality
holds iff $(n-1)/2$ of the $v_i$ are equal to some unit vector $v$
and another $(n-1)/2$ are equal to $-v$. This follows easily from
the proof above.

We mention further that replacing the condition $u\cdot v_i \ge 0$ by
$u\cdot v_i > 0$ for every $i\in [n]$ in
Theorem~\ref{semiespacio} does not imply $\|v_1 + v_2 + .
. . + v_n\| > 1.$ For instance when $\|.\|$ is the max norm and $v_1=\ldots =v_k=(-1,\eps)$
and $v_{k+1}=\ldots =v_{2k-1}=(1,\eps)$ and $\eps>0$ is small enough, $\|s\|=1$ although $u\cdot v_i>0$ for all $i$.

\noindent {\bf Remark 2.} Theorem~\ref{semiespacio} has no
analogue in dimension 3 and higher. For the example showing this
let $B$ be the Euclidean unit ball in $\R^3$, let $L$ be a plane
at distance $\eps$ from the origin with unit normal $u$, and let
$P_n$ be a regular $n$-gon inscribed in the circle $L\cap B$, with
vertices $v_1,\ldots,v_n$. It is clear that $u\cdot v_i>0$ for all
$i\in [n]$ but $\sum_1^nv_i=\eps nu$ whose norm is as small as you
wish. The parity of $n$ does not matter.

\noindent {\bf Remark 3.} The following is a direct consequence of
Theorem~\ref{semiespacio}: let $V=\{v_1,v_2, \ldots v_n\}$ be  a
set of  unit vectors in a normed plane. Then it is always possible
to choose numbers $\epsilon_1, \epsilon_2, \ldots, \epsilon_n$
from $\{1,-1\}$ such that for every subset $W\subset V$ of odd
size, we have that $\|\sum\limits_{v_i\in W} \epsilon_i v_i \|\geq
1.$

\section{Proof of Theorem~\ref{hellytype}}

We need some preparations before the proof. We start with a small piece from Euclidean plane geometry. Let
$a,b,c$ be distinct unit vectors in the Euclidean plane and define
$\triangle =\conv\{a,b,c\}$. It is well known $h =
a + b+ c$ is outside  $\triangle$ (indeed, outside
the unit circle) if the triangle is obtuse, and is inside
$\triangle$ if the triangle is acute. (We ignore right angle triangles here.)
This is equivalent to saying that $h \in \tri$ iff $0 \in \tri$ since $\tri$ is
acute or obtuse depending on whether $0 \in \tri$  or not.

Is this statement true for any norm in $\R^2$? As we
see from the following lemma the answer is yes.

\begin{lemma}\label{conv} Assume $a,b,c \in \partial B$ and set $\tri=\conv \{a,b,c\}$. Then $0 \in \tri$ if
and only if $h=a+b+c \in \tri$.
\end{lemma}

\noindent {\bf Proof.} If $0 \notin \tri$, then by separation
there is a vector $u$ such that $u\cdot a,u\cdot b,u\cdot c > 0$.
Theorem~\ref{semiespacio} with $V=\{a,b,c\}$ applies and shows
that $h\notin \intt B$. As $\intt \tri \subset \intt B$, $h \in
\tri$ implies $h \in \partial \tri$, say $h \in [a,c]$. Then
$a+b+c=ta+(1-t)c$ for some $t \in [0,1]$ and so
\[
\frac {1-t}2 a+\frac 12 b + \frac t2 c=0,
\]
a convex combination of $a,b,c$, showing that $0 \in \tri$. So indeed, $h \notin \tri.$

Assume next that $0 \in \tri$. Since $0$ is the center of the unit
ball, it must be contained in the medial triangle of $\tri$, that
is, $0=\alpha(\frac{b+c}{2})+\beta (\frac{a+c}{2})+\gamma(\frac{a+b}{2}),$ with
$\alpha,\beta,\gamma\in [0,1]$ and $\alpha+\beta+\gamma=1.$ We
have that
$$0+ \frac{\alpha}{2}\cdot a+\frac{\beta}{2}\cdot b+\frac{\gamma}{2}\cdot c=
\left(\alpha+\beta+\gamma\right)\left(\frac{a+b+c}{2}\right)=\frac{a+b+c}{2},$$

then $a+b+c=\alpha a+\beta b+\gamma c$, that is, $h=a+b+c\in \tri.$
\fin

\noindent {\bf Proof} of Theorem~\ref{hellytype}. We assume first
the extra condition that $V$ contains no antipodal pair of points.
For distinct $i,j,k \in [n]$, the vector $h=v_i+v_j+v_k$ is not in
$\intt B$. As all points of the triangle $\tri=\conv
\{v_i,v_j,v_k\}$ except for its vertices lie in $\intt B$, $h
\notin \tri$ either, unless $h=v_i$ say. But in this case $v_j$
and $v_k$ are antipodal.

So $h \notin \tri$ and then $0 \notin \tri$ follows from Lemma~\ref{conv}. Carath\'eodory's theorem
(see \cite{Car}) shows that $0 \notin \conv V$, too. By separation, there is a vector $u\neq 0$ with $u\cdot v>0$ for every $i \in [n]$.
Theorem~\ref{semiespacio} applies and gives $\|v_1+\ldots +v_n\|\ge 1$.

The general case goes by induction on $n$. The starting case $n=3$ is trivial. In the induction step $n-2 \to n$ (when $n\ge 5$)
$V=\{v_1,\ldots,v_n\}$ either satisfies the extra condition and we are done, or $V$ contains an antipodal pair, $v_{n-1},v_n$ say. By induction,
$\|v_1+\ldots +v_{n-2}\|\ge 1$, and the equality $\sum_1^n v_i=\sum_1^{n-2}v_i$ finishes the proof. \fin

\noindent {\bf Remark 4.} Theorem~\ref{hellytype} has no direct
analogue in $\R^3$. For instance if $V$ is the set of vertices of
a regular tetrahedron centered at the origin and inscribed in the
Euclidean unit ball, then every triple sum has (Euclidean) norm
$1$ yet the sum of the vectors is zero. A second example is when
is $u$ is a unit vector in $\R^3$ $v_1,v_2,v_3$ are the vertices
of a regular triangle in the plane orthogonal to $u$ and center at
$u$ and $v_{i+3}=v_i-2u$ ($i=1,2,3$), $V=\{v_1,\dots,v_6\}$ and
$B=\conv\{\pm v_1,\ldots,\pm v_6\}$. The sum of any three  vectors
from $V$ has norm at least one but $\sum_1^6 v_i=0$. The same
example works for Therorem~\ref{helly}, this time every 4-sum has
norm larger than one but $\sum _1^6v_i=0$ again.

\section{Preparations for the proof of Theorem~\ref{helly}}

We need a lemma about 6 vectors in the plane.

\begin{lemma}\label{main} Assume $z_1,\ldots,z_6 \in B$ and $\sum_1^6 z_i=0$. Then there are
distinct $i,j,k$ with $z_i+z_j+z_k \in B$.
\end{lemma}

\noindent {\bf Proof.} Assume for the time being that there are
two linearly independent vectors among the $z_i$. We will deal
with the remaining case soon. Define $D=\conv\{\pm z_1,\ldots,\pm
z_6\}$, $D$ is an $0$-symmetric convex polygon with at most 12
vertices. Clearly $Z=\{z_1,\ldots,z_6\} \subset D$ and $D \subset
B$. This implies that it suffices to prove Lemma~\ref{main} when
$B=D$.

Let $\ver D$ denote the set of vertices of $D$. We distinguish two cases:

{\bf Case 1.} When $|Z \cap \ver D|=2$. Then $D$ is a parallelogram with vertices $a,b,-a,-b$ where $a,b \in Z \cap \ver D$. As the assumptions and statement of the lemma are invariant under a non-degenerate linear transformation we may assume that $a=(1,1)$ and $b=(-1,1)$. This is in fact the case of the max norm. We need the following

\begin{claim} If the sum of real numbers $z_1,\dots,z_6$ is zero and all of them lie in $I=[-1,1]$, then there are at least 12 distinct triplets among them whose sum lies in $I$ as well.
\end{claim}

The proof is postponed to Section~\ref{claims}. We note first that
Claim~1 justifies our assumption about the existence of two
linearly independent vectors among the $z_i$. Indeed, if all the
$z_i$ are on a line through the origin, then they can be thought
of as real numbers. Claim~1 says then that there are three among
them with the required property (actually, 12 such triplets).

We show next how the Claim finishes Case 1. Both the first and the
second components of the $z_i$ satisfy the conditions of Claim 1.
So there are 12 triplets whose first components, and 12 further
triplets whose second component, sum to a number in $I$. As there
are 20 triplets altogether, there is a triplet whose first and
second components sum to a number in $I$, that is, there are
distinct $i,j,k$ with $z_i+z_j+z_k \in D$.

{\bf Case 2.} When $|Z \cap \ver D|\ge 3$. If there are $a,b,c \in
Z \cap \ver D$ such that $a+b+c \in D$, then we are done.
Otherwise Lemma~\ref{conv} (together with Carath\'eodory's
theorem) says that $0 \notin \conv\; (Z \cap \ver D)$. So we may
assume that every point of $Z \cap \ver D$ is in the open upper
halfplane. Let $a$ be the first and $b$ be the last vertex as we
walk around $\partial D$ in the upper halfplane in anticlockwise
direction. By a non-degenerate linear transformation we can
achieve $a=(1,1)$ and $b=(-1,1)$. Clearly, $[a,-b]$ and $[b,-a]$
lie on $\partial D$. Note that there is $c=(c_1,c_2) \in Z\cap
\ver D$, different from $a,b$ implying that $c_1 \in (-1,1)$ and
$c_2>1$.

\centerline{\psset{unit=1.7cm}
\begin{pspicture}(-2,-1.8)(2,2)
\pspolygon[linewidth=.02,fillcolor=green,fillstyle=solid](1,1)(-.3,1.4)(-1,1)(-1,-1)(.3,-1.4)(1,-1)
\psaxes{->}(0,0)(-2,-1.8)(2,1.8)
\psline[linestyle=dashed,linecolor=darkgray](0,1.4)(-.3,1.4)(-.3,0)
\rput(1.12,1.05){\small $a$} \rput(-1.12,1.05){\small $b$}
\rput(1.12,-1.05){\small $-b$} \rput(-1.14,-1.12){\small $-a$}
\rput(-.3,-.15){\small $c_1$} \rput(.15,1.4){\small $c_2$}
\end{pspicture}}
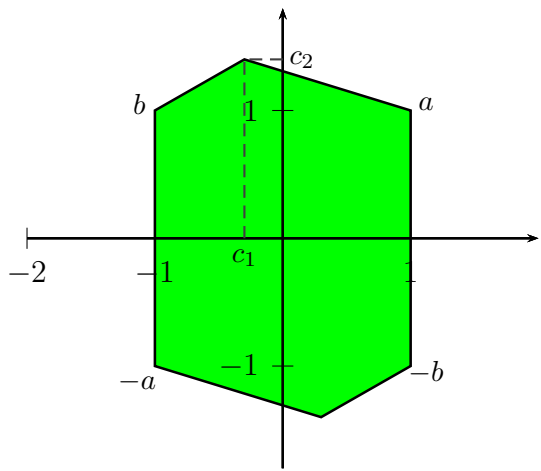
\centerline{\small {Figure 1.}}

For simpler writing let $u_1,u_2,u_3$ be the $z_i$ distinct from $a,b,c$. We are going to show that $a+b+u_i \in D$ for some $i$. Otherwise $u_i \notin D-a-b$ for all $i$. In other words, $u_1,u_2,u_3 \in D \setminus (D-a-b)$. It is easy to see that the second component of every vector in $D \setminus  (D-a-b)$ is larger than $-1$. Thus the second component of $u_1+u_2+u_3$ is larger than $-3$. The second component of $a+b+c$ is $2+c_2>3$. This contradicts the assumption $z_1+\ldots +z_6=0$.\fin

To close this section we prove Theorem~\ref{helly} in the case when $V$ does not contain two linearly independent vectors. In this case $V$ can be thought of as real numbers $x_1,\ldots,x_n$ with $x_1\ge \ldots \ge x_n$. By symmetry and scaling we may assume that $x_1=1\ge |x_n|$ and $B=[-1,1]$. There is nothing to prove if $x_n\ge 0$. Also, $x_1=-x_n$ is impossible since then $x_1+x_2+x_n=x_2 \in B$, contrary to the conditions. Thus $x_1+x_n>0$ and $x_{n-1}>0$ as otherwise $x_1+x_{n-1}+x_n \in B$. Consequently $x_1+\ldots +x_n\ge x_1+x_{n-1}+x_n>1$.

\section{Proof Theorem~\ref{helly}}

The result is trivially true for $n=3$. Next comes the case $n=5$: Set $z_i=v_i$, $i=1,2,3,4,5$ and $z_6=-(v_1+\ldots+v_5)$.
If $\|z_6\|\le 1$ were the case, then Lemma~\ref{main} implies that a 3-sum, $z_i+z_j+z_k$ say, lies in $B$. This contradicts the condition if $z_6$ is not present among $z_i,z_j,z_k$. But if it is, then the complementary 3-sum goes without $z_6$, and its norm equals $\|z_i+z_j+z_k\|\le 1$, a contradiction again.

Assume now that the theorem fails and let $V=\{v_1,\ldots,v_n\}$ be a counterexample with the smallest possible $n$ and let $B$ be the unit ball of the corresponding norm. Here $n\ge 7$ clearly and $V$ contains two linearly independent vectors. Define $v_0=\sum_1^nv_i$. Then $D=\conv \{\pm v_0,\pm v_1,\ldots,\pm v_n\}$ is an $0$-symmetric convex body (actually a convex polygon) that is the unit ball of a norm $\|.\|$. As $D \subset B$, $V$ is a counterexample with this norm. This means that $\|v_i\|\le 1$ for all $i=0,1,\ldots,n$ and all 3-sums have norm $>1$. From now on we keep this norm fixed and consider $V$ a counterexample with this norm.

We choose $\la<1$ but very close to $1$ so that $\la v_1,\dots,\la v_n$ is still a counterexample, this time with $\|\sum_1^n \la v_i\| < 1$.
By continuity there is an $\eps>0$ so that if $\|u_i-\la v_i\| < \eps$ for all $i \in [n]$, then $U=\{u_1,\dots,u_n\}$ is still a counterexample meaning that $\|u_i\|\ <1$ for all $i\in [n]$, $\| \si (S,U)\|>1$ for all $S \in {[n] \choose 3}$, and $\|\sum_1^n u_i\|<1$.
Here of course $\si (S,U)$ stands for $\sum_{i \in S}u_i$.

\begin{claim} One can choose $U$ so that for all $S \in {[n] \choose 3}$ and all $T \in {[n] \choose 3}\cup {[n]\choose 5}$, $\|\si (S,U)\|=\|\si (T,U)\|$ implies $S=T$.
\end{claim}

The technical proof is postponed to Section~\ref{claims}. Now we return to the proof by fixing $U$ as in the claim.

The numbers $\|\si (S,U)\|$ with $S \in {[n]\choose 3} \bigcup {[n] \choose 5}$ are all larger than one. Let $\mu>1$ be the smallest among them.
We claim that $\mu=\|\si (S,U)\|$ for some unique $S \in {[n] \choose 3}$. Indeed, if the minimal $S$ is a 5-tuple, $S=\{1,2,3,4,5\}$ say, then the five vectors $\mu^{-1}u_1,\ldots,\mu^{-1} u_5$ are all in $D$, all of their 3-sums are outside $D$ but their sum is in $D$, contradicting case $n=5$ of the theorem.

Consequently $\mu=\|\si (S,U)\|$ for a unique $S \in {[n] \choose 3}$. We assume w.l.o.g. that $S=\{1,2,3\}$. Choose $\nu <\mu^{-1}<1$ so that $\nu \|\si (T,U)\|>1$ for all $T \in {[n]\choose 3} \cup {[n] \choose 5}$ except for $T=S$ and $\nu\|\si(S,U)\|<1$. Set $w_0=\nu(u_1+u_2+u_3)$, $w_i=\nu u_i$ for $i>3$, and define $W=\{w_0,w_4,\dots,w_n\}$.

We show finally that $W$ is another counterexample with the norm $\|.\|$. This would contradict the minimality of $n$ as $|W|=n-2<n$ and so finish the proof.

It is clear that $W \subset D$ and $w_0+w_4+\ldots +w_n\in D$. All 3-sums of $W$ that do not contain $w_0$ are outside $D$ since such a 3-sum equals $\nu(u_i+u_j+u_k)$ with $4\le i<j<k$ which is outside $D$ by the definition of $\nu$. A 3-sum of the form $w_0+w_i+w_j$ for $4 \le i<j$ is equal to $\nu(w_1+w_2+w_3+w_i+w_j)$ which is again outside $D$ because of the definition of $\nu$.\fin

\section{Proofs of the Claims}\label{claims}

\noindent {\bf Proof} of Claim 1. Write $x_1,x_2,\ldots,x_s$ resp
$-y_1,\ldots,-y_t$ for the positive and non-positive elements of
our set $Z$ of real numbers, here $s+t=6$ and we assume w.l.o.g.
that $s\le t$. We assume further that $x_1\ge x_2 \ge \ldots \ge
x_s$ and $y_1\ge \ldots \ge y_t$. The case $s=0$ is trivial, and
so is case $s=1$: then all 3-sums of $Z$ lie in $I=[-1,1]$.

If $s=2$, then $x_1-y_i-y_j  \in I$ for all distinct $i,j$. Indeed, this is clear if $x_1\ge y_i+y_j$ since then $0\le x_1-y_i-y_j\le x_1 \le 1$. Assume next that $x_1<y_i+y_j$ and $y_i\ge y_j$ say, then $-1\le -y_i \le -y_i+(x_1-y_j)<0$ provided $x_1\ge y_j$. But case $x_1<y_j$ is impossible: then we'd have $x_1,x_2 < y_i,y_j$ and $x_1+x_2 < y_i+y_j$ so the sum of our six numbers cannot zero. Thus there are ${4 \choose 2}=6$ distinct 3-sums in $I$ and no two of them are complementary. The 6 complementary 3-sums lie in $I$, too.

Finally $s=3$. By symmetry we assume that $x_1\ge y_1$. If $x_1,x_2\ge y_1$, then $x_k-y_i-y_j \in I$ for $k=1,2$ and for all distinct $i,j$. This follows the same way as above. This is already 6 distinct 3-sums in $I$ (with no two complementary), giving 12 distinct 3-sums that lie in $I$.

So suppose $y_1>x_2$. Again $x_1-y_i-y_j \in I$ for all distinct $i,j$ and both $-y_1+x_1+x_2$ and $-y_1+x_1+x_3$ lie in $I$ as both are non-negative and each smaller than $x_1$. This is five distinct (and non-complementary) 3-sums. We only need to find one more.

The missing one is $-y_1+x_2-y_2$ if $x_1\ge y_1>x_2 \ge y_2$, and $x_1-y_2+x_2$ if $x_1\ge y_1\ge y_2 \ge x_2$. \fin

\bigskip
\noindent {\bf Proof} of Claim 2. Our unit ball $D$ is an
$0$-symmetric convex polygon with edge set $E$. For an edge
$e=[x,y]$ define $\ell_e$ as the (unique) linear function $\R^2
\to \R$ such that $\ell_e(x)=\ell_e(y)=1$. It follows that for all
$z \in \R^2$, $\|z\|= \min \{\ell_e(z): e \in E\}$.

Recall the definition of $[n] \choose k$ and $\si(S,V)$ from Section 1. We are going to choose the vectors $u_1,u_2,\ldots,u_n$ in this order where $u_i$ is in the $\eps$-neighbourhood $N_{\eps}(\la v_i)$ of $\la v_i$) ($i \in [n]$) so that the following holds. The sets $U_k=\{u_1,\ldots,u_k\}$ for $k\in [n]$ satisfy
\begin{enumerate}
\item[(1)] $\ell_e(\si (S,U_k))\ne \ell_f(\si (T,U_k))$ for all distinct $e,f \in E$ and all $s,t\in\{0,1,\ldots,5\}$ and all $S \in {[k] \choose s}$ and $T \in {[k] \choose t}$ with $S \ne T$,

\medskip
\item[(2)] $\ell_e(\si (S,U_k))\ne \ell_e(\si (T,U_k))$ for all $e \in E$, for all $s,t\in\{0,1,\ldots,5\}$ and all $S \in {[k] \choose s}$ and $T \in {[k] \choose t}$ with $S\ne T$.
\end{enumerate}

These conditions guarantee that in $U=U_n$ all $3$-sums have different norms and no 3-sum and 5-sum have the same norm.
This is the requirement in Claim 2.

The proof goes by induction. The first vector $u_1$ is chosen from
$N_{\eps}(\la v_1)$ so that $\ell_e(u_1)\ne 0$ for all $e \in E$. So the
forbidden region for $u_1$ is the union of finitely many lines, and
consequently there is a suitable $u_1$. Assume $U_k$ has been constructed satisfying conditions (1) and (2) and $k\ge 1$.

We start with condition (1). For a fixed pair $e,f \in E$ ($e\ne f$), and for a fixed $S \in {[k+1] \choose s}$ and fixed $T \in {[k+1] \choose t}$, (1) says something for $u_{k+1} \in N_{\eps}(\la v_{k+1})$ only if $k+1 \in S\cup T$, otherwise it is satisfied by the induction hypothesis. If $k+1$ only appears in $S$, (resp. in $T$), then (1) says that $\ell_e(u_{k+1}) \ne \al$ (and $\ell_f(u_{k+1})\ne \al$) for a particular value of $\al$ depending only on $e,f,S,T$. So the forbidden region is a line $L=L(e,f,S,T)$. When $k+1 \in S \cap T$ then the condition is $\ell_e(u_{k+1})-\ell_f(u_{k+1})\ne \al$. So the forbidden region is a line again as $\ell_e-\ell_f$ is a non-identically zero linear function.

Checking condition (2) is similar. For a fixed $e \in E$, and for fixed $S \in {[n] \choose s}$ and $T \in {[n] \choose t}$, condition (2) says something for $u_{k+1}$ only if again $k+1 \in S \cup T$, otherwise it is satisfied by the induction hypothesis. If $k+1 \in S \cap T$, then condition (2) says that $\ell_e(\si (S\setminus\{k+1\},U_{k+1}))\ne \ell_e(\si (T\setminus \{k+1\},U_{k+1})$. This follows from the induction hypothesis. Finally, if $k+1$ is in $S\setminus T$, condition (2) says that $\ell_e(u_{k+1}) \ne \al$ with a particular value of $\al$ depending only on $e,f,S,T$. So the forbidden region is a line, again. The same applies when $k+1\in T \setminus S$.

As there are finitely many such forbidden lines for $u_{k+1}$, the Lebesgue measure of the forbidden region is zero. Thus almost all choices of $u_{k+1}$ avoid the forbidden region. \fin

\section{Characterization of central symmetry}

Theorem~\ref{hellytype} is about a norm whose unit ball is an
$0$-symmetric convex body $B$. In the particular case $n=3$ it says that if $a,b,c$
are unit vectors and their convex hull is separated from $0$, then
their sum has norm at least 1. The next theorem is a kind of
converse.

\begin{theorem} Let $K\in \mathbb R^2$ be a convex body with $0 \in {\emph\intt} K$. Then $K$ is
centrally symmetric with center at $0$ under either one of the following conditions.
\begin{enumerate}
\item[\rm(i)] For any three distinct vectors $a,b,c\in
\partial K$ contained in a closed halfplane whose bounding line
goes through $0$, the vector $a+ b +c \notin {\emph\intt} K$.
\item[\rm(ii)] For any three distinct vectors $a,b,c \in
\partial  K$ with $0 \in \emph{int conv\;}\{a,b,c\}$, the vector
$a + b + c \in {\emph\intt} K$
\end{enumerate}
\end{theorem}

\noindent {\bf Proof} of (i). Suppose on the contrary that $K$ is
not centrally symmetric. Then we can choose a chord $ac$ (of $K$)
containing $0$ with $a+c\neq 0$. Further let $b$ be a vector on
$\partial K,$ very close to $a$, and let $bw$ be the chord which
is parallel to $ac$. It is very easy to see that $h=a+b+c\in
\text{relint}(bw)$ if $b$ is close enough to $a$. This implies
that $h\in \text{int}K$, a contradiction. \fin

\centerline{\psset{unit=.8cm}
\begin{pspicture}(0,.5)(6,6.4)
\definecolor{rosita}{cmy}{0.1,.3,.1}%
\pscurve[linewidth=.035,fillcolor=rosita,fillstyle=solid](6,3)(5.4,5.5)(2.4,5.5)(1,3)(2.5,1)(5,1)(6,3)
\psline(1.05,3.6)(5.98,3.6)
\psline[linecolor=blue,linewidth=.035]{->}(3,3)(6,3)
\psline[linecolor=blue,linewidth=.035]{->}(3,3)(1,3)
\psline[linecolor=blue,linewidth=.035]{->}(3,3)(1.05,3.6)
\psline[linestyle=dashed,linecolor=gray](4,3)(2.05,3.6)
\psline[linecolor=red,linewidth=.035]{->}(3,3)(2.05,3.6)
\psdots[dotsize=3pt](3,3) \rput(3,2.7){\small
$0$}\rput(.7,3){\small $a$}\rput(6.3,3){\small
$c$}\rput(.8,3.7){\small $b$}\rput(2,3.85){\small $h$}
\rput(6.2,3.7){\small $w$}
\end{pspicture}}
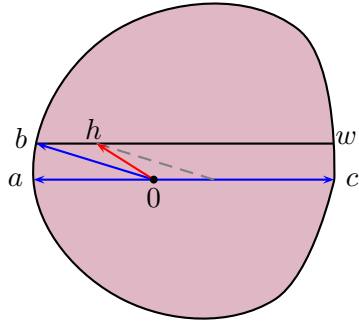
\centerline{\small {Figure 2.}}

\noindent {\bf Proof} of (ii). Again, let $ac$ be a chord of $K$
containing $0$ such that $a+c \ne 0$ and further, let $b$ be the
point on $\partial K$ where the tangent line $\ell$ at $b$ to $K$
is parallel to $ac$. We choose $a$ and $c$ so that this $b$ is a
single point (on either side of the chord $ac$). This is clearly
possible.

This way $h=a+b+c\in \ell$ and consequently $h$ is outside $K$.
Now, replace $a$ resp. $c$, by $a_1$ and $c_1$ very close to $a$
and $c$ so that the chord $a_1c_1$ is parallel to $\ell$ and so
that the line through $a$ and $c$ separates $b$ and $a_1,c_1.$ In
this case $0 \in \text{int}(\text{conv}\{a_1,b,c_1\})$. Since the
norm of the sum of vectors is a continuous function, we have that
$h_1=a_1+b+c_1$ is not in $\text{int} K$ provided the line
through $a_1c_1$ is close enough to the chord $ac$. \fin

\centerline{\psset{unit=.88cm}
\begin{pspicture}(0,.5)(6,6.8)
\definecolor{azulado}{cmy}{0.5,0,.1}
\pscurve[linewidth=.035,fillstyle=solid,fillcolor=azulado](3,6)(2,5.6)(1,3)(2.5,1)(4,1)(5.5,3)(5,5)(3,6)
\psline(1,3)(5.5,3)\psline[linecolor=magenta]{->}(2.5,3)(3,6)
\psline[linecolor=magenta]{->}(2.5,3)(1.05,2.8)\psline[linecolor=magenta]{->}(2.5,3)(5.45,2.8)
\psline[linecolor=red,linewidth=.035]{->}(2.5,3)(4.02,5.87)
\psline[linestyle=dashed,linecolor=blue](3,6)(1.05,2.8)(5.45,2.8)(3,6)
\psline(1,6)(6,6)\psline[linestyle=dashed,linecolor=darkgray](4,3)(4.5,6)
\psdots[dotsize=3pt](2.5,3)(3,6)(1.05,2.8)(5.45,2.8)(4.5,6)
\rput(4.27,5.85){\small $h_1$} \rput(4.5,6.3){\small $h$}
\rput(2.3,3.25){\small $0$}\rput(.7,3.1){\small
$a$}\rput(5.8,3.1){\small $c$}\rput(.75,2.7){\small
$a_1$}\rput(5.8,2.7){\small $c_1$}\rput(3,6.25){\small
$b$}\rput(1.5,6.2){\small $\ell$}
\end{pspicture}}
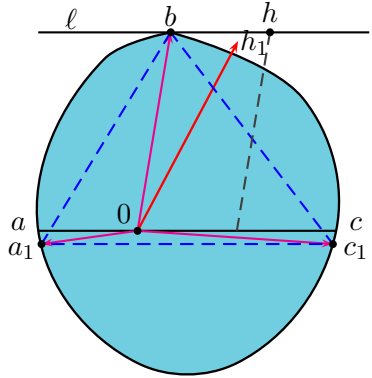
\centerline{\small {Figure 3.}}

\bigskip
{\bf Acknowledgements.} The authors are indebted to Viktor Grinberg for comments and discussions,
and in particular for the question that led from Theorem~\ref{hellytype} to Theorem~\ref{helly}.
The authors acknowledge the generous support of the Hungarian-Mexican
Intergovernmental S\&T Cooperation Programme T\'ET\_10-1-2011-0471 and NIH
B330/479/11 ``Discrete and Convex Geometry''. Research of the first author was partially
supported by ERC Advanced Research Grant no 267165 (DISCONV), and
by Hungarian National Research Grant K 83767.

\vspace{1cm} {\sc Imre B\'ar\'any}
\\
  {\footnotesize R\'enyi Institute of Mathematics}\\
  {\footnotesize Hungarian Academy of Sciences}\\
  {\footnotesize POBox 127, 1364 Budapest, Hungary}\\
  {\footnotesize e-mail: {\tt barany@renyi.hu}\\
  {\footnotesize and}\\
  {\footnotesize Department of Mathematics}\\
  {\footnotesize University College London}\\
  {\footnotesize Gower Street, London WC1E 6BT}\\
  {\footnotesize England}\vspace{.5cm}
\\
{\sc Jes\'us Jer\'onimo-Castro}
\\
  {\footnotesize Facultad de Ingenier\'ia,}\\
  {\footnotesize Universidad Aut\'onoma de Quer\'etaro}\\
  {\footnotesize Cerro de las Campanas s/n, C.P. 76010}\\
{\footnotesize Quer\'etaro, M\'exico}\\
   {\footnotesize e-mail: {\tt jesusjero@hotmail.com}}

\end{document}